\documentclass[12pt,leqno]{article}

\usepackage[cp1250]{inputenc}
\usepackage[OT4]{fontenc}
\usepackage{amsthm}
\usepackage{amsfonts}
\usepackage{amsmath}
\usepackage{amssymb}

\newcommand{\R}{\mathbb{R}}
\newcommand{\C}{\mathbb{C}}

\newcommand{\bo}{{\bf 0}}

\newcommand{\Aa}{\mathcal{A}}
\newcommand{\K}{\mathbb{K}}

\newcommand{\Opis}{{\cal O}}

\newtheorem{theorem} {Theorem}[section]
\newtheorem{theorem*}{Theorem}
\newtheorem{prop*} {Proposition}
\newtheorem{lemma*}{Lemma}
\newtheorem{lemma}[theorem]{Lemma}
\newtheorem{cor}[theorem]{Corollary}
\newtheorem{cor*}{Corollary}
\newtheorem{prop}[theorem] {Proposition}

\newtheorem{definition*}{Definition}

\theoremstyle{definition}

\begin{document} 
 

\begin{center}
{\Large Counting indices of critical points of rank two\\ of  polynomial
selfmaps of $\R^4$}\\[1em]
by  Zbigniew Szafraniec\\[1em]
\end{center}
{\bf Abstract} For a generic $f\in C^\infty(\R^4,\R^4)$ there is a discrete set  $\Sigma^2(Df)$ of critical
points of rank two, and there is an integer index $I_p(Df)$ associated to any $p\in\Sigma^2(Df)$.
We  show how to compute $\sum I_p(Df)$, $p\in\Sigma^2(Df)$, in the case where $f$
is a polynomial mapping.\\[1em]
\section{Introduction}
Assume that $f:M\rightarrow N$ is a $C^1$-mapping between 4-dimensional oriented manifolds. We shall denote by
$\Sigma^2(Df)$ the set of those critical points of $f$ where the derivative $Df$ is of rank two.
According to the Thom transversality theorem, critical points of a generic $f$ have rank $\geq 2$,
and the set $\Sigma^2(Df)$ is discrete.

Several authors observed that 
one may associate an index $I_p(Df)$ to each $p\in\Sigma^2(Df)$. If $M$ is a closed manifold then the algebraic sum
of indices $\#\Sigma^2(Df)=\sum_p I_p(Df)$, where $p\in\Sigma^2(Df)$, is a natural  invariant associated to $f$
(see  \cite{arnoldetal, ekholmtakase, saekisakuma, saji, stingley, takase}, and \cite{pragaczweber} for the comples case).
R.~Stingley \cite{stingley} proved that $\#\Sigma^2(Df)=\deg(f)\cdot p_N-p_M$, where $\deg(f)$ is the topological degree
of $f$, and $p_N$ (resp. $p_M$) is the Pontryagin number of $N$ (resp. $M$). This result demonstrates that there
is a non-trivial linear relation between those two natural invariants, i.e. $\#\Sigma^2(Df)$ and $\deg(f)$, associated to $f$.

In this paper we show how to compute $\#\Sigma^2(Df)$ in the case where $f:\R^4\rightarrow\R^4$ is such a polynomial mapping
that the family of critical points of rank two of its complexification is finite.
Our approach does not exclude the case where points in $\Sigma^2(Df)$ are not
umbilic points.

Papers \cite{ekholmtakase, saji, stingley} offer  methods of
computing the index $I_p(Df)$ which require the germ $f:(\R^4,p)\rightarrow (\R^4,f(p))$ to be written in a special form.
However, if $f$ is a polynomial then we usually are not  able to find explicitely
points in $\Sigma^2(Df)$, so we cannot adopt these techniques in our case.

Let $L$ denote the space of $4\times 4$-matrices, and let  $\Sigma$ denote the connected oriented $12$-dimensional submanifold of $L$
 consisting of matrices of rank two. As $Df:\R^4\rightarrow L$, one may define the index  $I_p(Df)$
as the intersection number of the mapping $Df$ with $\Sigma$ at $Df(p)\in\Sigma$. This is why
in Sections 2,3 we investigate a more general case of mappings $m:\R^4\rightarrow L$, and for $p$ isolated in $ m^{-1}(\Sigma)$ we introduce
the index $I_p(m)$.

In Section 4 we explain how to verify whether $m_{\R}^{-1}(\Sigma)$ is finite  in the case where  $m_{\R}:\R^4\rightarrow L$  is a polynomial mapping. Then we construct
a quadratic form whose signature equals $\#\Sigma^2(m_{\R})=\sum I_p(m_{\R})$, $p\in m_{\R}^{-1}(\Sigma)$.
In the end of this section we present examples
 which were calculated with the help
of {\sc Singular} \cite{greueletal}. 
We also give simple examples which demonstrate that there is only a trivial
linear relation between the two natural integer invariants, i.e. $\#\Sigma^2(Df)$ and $\deg(f)$,  which one may associate to a  polynomial $f$.

It is proper to add that \cite{bobowikszafraniec, krzyzanowskaszafraniec} present methods for counting the signed  cusp
or swallowtail singularities for polynomial selfmaps of $\R^n$, where $n=2,3$.

\section{Permutations of rows or columns}

In this section we show that  a $3$-cyclic permutation of either rows or columns of a $4\times 4$-matrix may be prolonged to
an isotopy preserving the rank of matrices.

Let $W$ be a vector space, and let $w_1,\ldots,w_4\in W$. Put
$w_1(t)=(1-t)w_1+tw_2$, $w_2(t)=(1-t)w_2+tw_3,$
$w_3(t)=(1-t)w_3+tw_1$, $w_4(t)=w_4.$
This way there is given an isotopy
$$\prod_{i=1}^4 W\times [0,1]\ni ((w_i)_{i=1}^4,t)\mapsto (w_i(t))_{i=1}^4\in \prod_{i=1}^4 W\, .$$
In particular $\left(  w_i(0)  \right)_{i=1}^4=(w_1,w_2,w_3,w_4)$, and $\left(  w_i(1)    \right)_{i=1}^4=(w_2,w_3,w_1,w_4)$
is a $3$-cyclic permutation of the first three vectors.

By the rank of a sequence of vectors we shall denote the dimension of the subspace spanned
by these vectors.

\begin{lemma}\label{per1}
$\operatorname{rank}\left(  w_i \right)_{i=1}^4=4$ if and only if $\operatorname{rank}\left(  w_i(t) \right)_{i=1}^4=4$ 
for all $0\leq t\leq 1$.
\end{lemma}
\noindent{\em Proof.} As $(1-t)^3+t^3>0$ for all $0\leq t\leq 1$, then the exterior product
$$w_1(t)\wedge w_2(t)\wedge w_3(t)\wedge w_4(t)=(1-t)^3w_1\wedge w_2\wedge w_3\wedge w_4+t^3w_2\wedge w_3\wedge w_1\wedge w_4$$
$$= ((1-t)^3+t^3)w_1\wedge w_2\wedge w_3\wedge w_4$$
does not vanish if and only if $w_1\wedge w_2\wedge w_3\wedge w_4\neq 0$. $\Box$

\begin{lemma}\label{per2}
$\operatorname{rank}\left(  w_i \right)_{i=1}^4\geq 3$ if and only if $\operatorname{rank}\left(  w_i(t) \right)_{i=1}^4\geq 3$ 
for all $0\leq t\leq 1$.
\end{lemma}
\noindent{\em Proof.} It is enough to prove $(\Rightarrow)$.
If $w_1,w_2, w_3$ are linearly independent then $w_1\wedge w_2\wedge w_3\neq 0$, and so
$$w_1(t)\wedge w_2(t)\wedge w_3(t)=((1-t)^3+t^3)w_1\wedge w_2\wedge w_3\neq 0$$
for all $0\leq t\leq 1$.

If $w_1,w_2,w_3$ are linearly dependent, then at least one of the products
$w_1\wedge w_2\wedge w_4$, $w_1\wedge w_3\wedge w_4$, $w_2\wedge w_3\wedge w_4$ is non-zero.
We have
$$
\left[ \begin{array}{c} w_1(t)\wedge w_2(t)\wedge w_4\\ w_1(t)\wedge w_3(t)\wedge w_4\\ w_2(t)\wedge w_3(t)\wedge w_4                        \end{array}        \right]\ 
=\ \left[ \begin{array}{ccc} (1-t)^2&t(1-t)&t^2\\ -t^2&(1-t)^2&t(1-t)\\ -t(1-t)&-t^2& (1-t)^2                        \end{array}     \right]
\left[ \begin{array}{c} w_1\wedge w_2\wedge w_4 \\ w_1\wedge w_3\wedge w_4\\ w_2\wedge w_3\wedge w_4                 \end{array}       \right].
$$
The determinant of the above matrix equals $((1-t)^3+t^3)^2$, so it does not vanish in $[0,1]$. Hence, for any
$0\leq t\leq 1$, at least one of the products $w_i(t)\wedge w_j(t)\wedge w_4$ does not vanish, and then
$\operatorname{rank}(w_i(t))_{i=1}^4\geq 3$. $\Box$

\begin{lemma}\label{per3}
$\operatorname{rank}\left(  w_i \right)_{i=1}^4\geq 2$ if and only if $\operatorname{rank}\left(  w_i(t) \right)_{i=1}^4\geq 2$ 
for all $0\leq t\leq 1$.
\end{lemma}
\noindent{\em Proof.} Assume first that $\operatorname{rank}(w_1,w_2,w_3)\geq 2$. Then at least one of the products
$w_1\wedge w_2$, $w_1\wedge w_3$, $w_2\wedge w_3$ is non-zero. We have
$$
\left[ \begin{array}{c} w_1(t)\wedge w_2(t)\\ w_1(t)\wedge w_3(t)\\ w_2(t)\wedge w_3(t)                       \end{array}        \right]\ 
=\ \left[ \begin{array}{ccc} (1-t)^2&t(1-t)&t^2\\ -t^2&(1-t)^2&t(1-t)\\ -t(1-t)&-t^2& (1-t)^2                        \end{array}     \right]
\left[ \begin{array}{c} w_1\wedge w_2 \\ w_1\wedge w_3\\ w_2\wedge w_3\end{array}\right] .$$
By the same arguments as in Lemma \ref{per2}, $\operatorname{rank}(w_1(t),w_2(t),w_3(t))\geq 2$ for all $0\leq t\leq 1$.

If $\operatorname{rank}(w_1,w_2,w_3)\leq 1$, then at least one of the products
$w_1\wedge w_4$, $w_2\wedge w_4$, $w_3\wedge w_4$ is non-zero. We have

$$
\left[ \begin{array}{c} w_1(t)\wedge w_4\\ w_2(t)\wedge w_4\\ w_3(t)\wedge w_4                     \end{array}        \right]\ 
=\ \left[ \begin{array}{ccc} 1-t&t&0\\ 0&1-t&t\\ t& 0& 1-t                       \end{array}     \right]
\left[ \begin{array}{c} w_1\wedge w_4 \\ w_2\wedge w_4\\ w_3\wedge w_4\end{array}\right] .$$
The determinant of the above matrix equals $(1-t)^3+t^3$. Hence for any $0\leq t\leq 1$ at least one of the products
$w_i(t)\wedge w_4$ does not vanish,
and then  $\operatorname{rank}(w_i(t))_{i=1}^4\geq 2$. $\Box$

\begin{lemma}\label{par4}
$\operatorname{rank}\left(  w_i \right)_{i=1}^4\geq 1$ if and only if $\operatorname{rank}\left(  w_i(t) \right)_{i=1}^4\geq 1$ 
for all $0\leq t\leq 1$.
\end{lemma}
\noindent{\em Proof.} If $w_4\neq 0$ then the assertion is obvious. Suppose that at least one $w_i\neq 0$, where $1\leq i\leq 3$.
Then
$$
\left[ \begin{array}{c} w_1(t)\\ w_2(t)\\ w_3(t)                    \end{array}        \right]\ 
=\ \left[ \begin{array}{ccc} 1-t&t&0\\ 0&1-t&t\\ t& 0& 1-t                       \end{array}     \right]
\left[ \begin{array}{c} w_1\\ w_2\\ w_3\end{array}\right] ,$$
and we may apply the same arguments as in the previous lemma. $\Box$

\begin{cor}\label{par5}
Obviously, if $\operatorname{rank}\left( w_i(t)\right)_{i=1}^4\geq k$ for at least one $t$,
then $\operatorname{rank}\left( w_i \right)_{i=1}^4\geq k$. Hence
$\operatorname{rank}\left(  w_i \right)_{i=1}^4=k$ if and only if $\operatorname{rank}\left(  w_i(t) \right)_{i=1}^4=k$ 
for each $0\leq t\leq 1$.
\end{cor}

Let $L$ denote the linear space of $4\times 4$-matrices
$$M\ =\ \left[\begin{array}{cccc} 
a_{11}&a_{12}&b_{11}&b_{12} \\
a_{21}&a_{22}&b_{21}&b_{22} \\

c_{11}&c_{12}&d_{11}&d_{12} \\
c_{21}&c_{22}&d_{21}&d_{22} 

\end{array}\right]\ =\ 
\left[\begin{array}{cc} A& B\\ C& D

\end{array}\right] $$
with real coordinates $(a_{11},a_{12},a_{21},a_{22},b_{11},\ldots ,c_{22},d_{11}, d_{12},d_{21},d_{22})$.

There is a natural isomorphism $L\simeq \prod_{i=1}^4\R^4=\{ \left(w_i\right)_{i=1}^4\ |\ w_i\in\R^4\}$, where the sequence $(w_i)_{i=1}^4$ denotes
either rows or columns of $M$. We get immediately
\begin{cor}\label{per6}
Let $\tau:L\rightarrow L$ be a composition of  finite sequence of $3$-cyclic permutations of either rows or columns.
(It is worth to notice that  a permutation consisting of two disjoint transpositions is a composition of a finite sequence of
$3$-cyclic permutations.)
Then there exists an isotopy $T:L\times [0,1]\rightarrow L$ such that $T(M,0)\equiv M$, $T(M,1)=\tau(M)$, and
$\operatorname{rank}\, T(M,t)\equiv \operatorname{rank}\, M$.

\end{cor}

\section{Index of a critical point}

In this section we  investigate the intersection number of a germ $(\R^4,p)\rightarrow L$ with the submanifold consisting of
matrices of rank two, and then  we introduce the index of a critical point of rank two.

Denote by $U_\pm$ the open connected subsets of $L$ consisting of matrices with $\pm\det(A)> 0$, and let $U=U_+\cup U_-$.
Put  $\Sigma=\{M\in L\ |\ \operatorname{rank}(M)=2\}$. By \cite[Proposition 2.5.3]{golub}, a matrix $M\in U$ has rank two
if and only if $D-CA^{-1}B=0$.
Then $\Sigma\cap U$ is the graph of the function $D=CA^{-1}B$, so that $\dim \Sigma\cap U=12$,
$\operatorname{codim} \Sigma\cap U=4$. There is also  given the natural orientation of the normal bundle over $\Sigma\cap U_+$ induced
by the orientation of coordinates $(d_{11},d_{12},d_{21},d_{22})$ in the source of this function. 

Let $M_{ij}$, where $3\leq i,j\leq 4$, denote the minor obtained by removing the $i$-th row and the $j$-th column 
from  $M$.
One may check that
$$\det(A)\cdot(D-CA^{-1}B)\ =\ 
\left[   \begin{array} {cc} M_{44}&M_{43}\\ M_{34}&M_{33}       \end{array}     \right],$$
so that $\Sigma\cap U=\bigcap M_{ij}^{-1}(0)\cap U$.
Moreover
$$\frac{\partial(M_{44},M_{43},M_{34},M_{33})  }{\partial ( d_{11},d_{12},d_{21},d_{22}     )   }\ =\ (\det(A))^4$$
is positive on $U$. In the further part of this paper we will need

\begin{lemma}\label{zamianka}
Let $M'$ denote the matrix obtained by interchanging the first row of $M\in U_+$ with the second one, and the third row with
the fourth one, so that $M'\in U_-$.  Let $M_{ij}'$ denote the minor obtained by removing the $i$-th row and the $j$-th column from  $M'$.
Then
$$
\left[ \begin{array}{c}M_{44}'\\ M_{43}'\\ M_{34}'\\ M_{33}'                 \end{array}   \right]\ =\ 
\left[ \begin{array}{rrrr}  0&0&-1&0\\ 0&0&0&-1\\ -1&0&0&0\\ 0&-1&0&0     \end{array}    \right]\,
\left[ \begin{array}{c}M_{44}\\ M_{43}\\ M_{34}\\ M_{33}                 \end{array}   \right]\ ,
$$
and the determinant of the above $4\times 4$-matrix equals $+1$. $\Box$
\end{lemma}

\begin{lemma}\label{minory_1}
Let $\cal{R}$ denote the localization of the ring of polynomials on $L$ by the powers of $\det(A)$. Then the ideal in $\cal{R}$ generated by all
$3\times 3$-minors of $M$ equals by the one generated by $M_{ij}$, where $3\leq i,j\leq 4$.

\end{lemma}
\noindent{\em Proof.} Applying elementary operations on rows and collumns, one may transform the matrix M to the form
$$
M''\ =\ \left[ \begin{array}{cccc}  a_{11}&a_{12}&0&0\\ a_{21}&a_{22}&0&0\\
0&0&M_{44}/\det(A)&M_{43}/\det(A)  \\0&0&M_{34}/\det(A)&M_{33}/\det(A)                  \end{array}      \right] \ .
$$
By the Cauchy-Binet formula, the ideal in $\cal{R}$ generated by all $3\times 3$-minors of $M$ is equal to the one generated by all $3\times 3$-minors
of $M''$, i.e. by all $M_{ij}$ and all $a_{ij}\cdot( M_{33}M_{44}-M_{34}M_{43} )/\det(A)^2$.
As $\det(A)$ is invertible in $\cal{R}$, the last ideal is generated by all $M_{ij}$. $\Box$

It is well-known that $\Sigma$ is a connected submanifold of $L$ of codimension $4$. According to \cite[Proposition 4.1]{ando},
 the manifold $\Sigma$, as well as its normal bundle,
is orientable.  Let fix the global orientation of of the normal bundle over $\Sigma$ which coincides with the orientation of 
the normal bundle over $\Sigma\cap U_+$ introduced before.\\[1em]
{\bf Definition.} Let $V$ be an open neighbourhood of $p\in\R^4$, and let $m:V\rightarrow L$ be a continuous mapping such that $p$ is isolated in
 $m^{-1}(\Sigma)$.
We define the index $I_p(m)$ as  the intersection number of $m$ with $\Sigma$ at $m(p)$.
In particular, if $m(p)\in \Sigma\cap U_+$ then
$I_p(m)$ is the local topological degree of the mapping 
$(\R^4,p)\ni x\mapsto H(x)=(h_{44}(x),h_{43}(x),h_{34}(x),h_{33}(x))\in (\R^4,\bo)$,
where $h_{ij}(x)=M_{ij}(m(x))$. 
By Corollary \ref{per6} and Lemma \ref{zamianka}, if $m(p)\in\Sigma\cap U_-$, then $I_p(m)$ equals the local topological degree
of the same mapping. If $p=\bo$ is the origin in $\R^4$, we shall denote its index by $I(m)$.\\[1em]

In the remainder of this section and in the next one we shall assume that $p=\bo$ is isolated in $m^{-1}(\Sigma)$.

Of course, if there is a continuous family of mappings $m_t:V\rightarrow L$, where $t\in[0,1]$, such that
$m_t(x)\in \Sigma$ if and only if $x=\bo$,
then $I(m_0)=I(m_1)$.
\begin{prop}\label{per7}
Assume that $\tau:L\rightarrow L$ is a composition of a finite sequence of $3$-cyclic permutations of rows
or columns such that $\tau( m(\bo))\in \Sigma\cap U$.
Then $I(m)=I(\tau\circ m)$.
\end{prop}
\noindent{\em Proof.}  We may assume that $m^{-1}(\Sigma)\cap V=\{\bo\}$.
Take such an isotopy $T:L\times [0,1]\rightarrow L$ as in Corollary \ref{per6}. Put $m_t(x)=T(m(x),t)$.  
Then $m_0=T(m,0)=m$, $m_1=T(m,1)=\tau\circ m$, so that $m_1(\bo)\in \Sigma\cap U$. Moreover
$m_t(x)\in \Sigma$ if and only if $x=\bo$.
Hence $I(m)=I(m_1)=I(\tau\circ m)$. $\Box$\\[1em]
{\bf Example.} Let
$$m(x,y,z,w)\ =\ \left[  \begin{array}{cccc}    x& y& z& 0\\
z^3& w& 0& 0\\ 0&0&1-x&y\\ w& 0&0 & 1-z                             \end{array}      \right].$$
Then $m(\bo)\in \Sigma\setminus U$. After applying a finite sequence $\tau$ of $3$-cyclic permutations
of rows or columns we get
$$\tau\circ m(x,y,z,w)\ =\ \left[  \begin{array}{cccc}    1-x& y& 0&0\\ 0&1-z&w& 0\\
z& 0&x& y\\ 0& 0& z^3& w        \end{array}      \right],$$
so that $\tau(m(\bo))\in \Sigma\cap U_+$. Then

$$\ \ \ \ \ \ \ \ \ h_{44}\ =\ \left[ \begin{array}{ccc}1-x&y&0\\0&1-z&w\\z&0&x               \end{array}    \right]\ =\ (1-x)(1-z)x+yzw,$$

$$h_{43}\ =\ \left[ \begin{array}{ccc} 1-x&y&0\\0&1-z&0\\z&0&y              \end{array}    \right]\ =\ (1-x)(1-z)y,$$

$$\ \ h_{34}\ =\ \left[ \begin{array}{ccc} 1-x&y&0\\0&1-z&w\\0&0&z^3              \end{array}    \right]\ =\ (1-x)(1-z)z^3,$$

$$h_{33}\ =\ \left[ \begin{array}{ccc}1-x&y&0\\0&1-z&0\\0&0&w               \end{array}    \right]\ =\ (1-x)(1-z)w.$$
Since the local topological degree $\deg_{\bo} ( H)=+1$, we have $I(m)=+1$.

\section{Polynomial mappings}
\noindent Let $\K$ denote either $\R$ or $\C$. For $p\in\K^4$, let $\Opis_{\K , p }$ denote the ring of germs
at $p$ of analytic functions $(\K^4,p)\rightarrow\K$, and let $\K[x]=\K[x_1,\ldots,x_4]$.
Let $L_{\C}$ denote the space of $4\times 4$-matrices having complex entries, and let $\Sigma_{\C}=\{M\in L_{\C}\ |\ \operatorname{rank}\,(M)=2\}$.

\begin{lemma}\label{przestawienia}
Assume that matrices $M_1,\ldots,M_s$ belong to $\Sigma_{\C}$. There exists an open dense
subset $\Delta\subset L\times L$ such that for any $(L_1,L_2)\in\Delta$ and each $1\leq i\leq s$,
the leading principal minor of $L_1\cdot M_i\cdot L_2$ of order $2$ does not vanish.
\end{lemma}
\noindent{\em Proof.} For $M\in L_{\C}$, let $a(M)$ denote its leading principal minor
of order 2. Then $a:L_{\C}\rightarrow \C$ is a polynomial mapping.

Each matrix $M_i$ is of rank 2, so that at least one of its $2\times 2$-minors does not vanish.
There exists $(L_1^i,L_2^i)\in L\times L$, which represents the appropriate interchange
of rows and columns, such that $a(L_1^i\cdot M_i\cdot L_2^i)\neq 0$. Then
$\Delta_i=\{(L_1,L_2)\in L\times L\ |\ a(L_1\cdot M_i\cdot L_2)\neq 0\}$ is a non-empty
Zarisky open subset of $L\times L$. The set $\Delta=\bigcap_1^s \Delta_i$ satisfies the assertion. $\Box$\\[1em]

Let
$$m_{\R}(x)\ =\ \left[   \begin{array}{cc} A(x)& B(x)\\ C(x)& D(x)  \end{array}       \right]\ :\ \R^4\rightarrow L$$
be a polynomial mapping, and let $m_{\C}:\C^4\rightarrow L_{\C}$ be its complexification. Let $H_{\K}=(h_{44},h_{43},h_{34},h_{33}):\K^4\rightarrow\K^4$, 
where $h_{ij}(x)=M_{ij}(m_{\K}(x))$.

Denote by $S_{\K}$ the ideal in $\K[x]$ generated by all $3\times 3$-minors of $m_{\K}(x)$. Set $\Aa_{\K}=\K [x]/S_{\K}$.
For $p\in V(S_{\K})$, denote $\Aa_{\K,p}=\Opis_{\K,p}/S_{\K}$.
From now on we shall assume that that $\dim_{\R}\Aa_{\R}=\dim_{\C}\Aa_{\C}<\infty$.
Then the set $V(S_{\K})$ of zeros of $S_{\K}$ in $\K^4$, which equals
$\{ p\in \K^4\ |\ \operatorname{rank}(m_{\K}(p))\leq 2\}$, is finite. 

Let $P_{\K}\subset\K[x]$
denote the ideal generated by all $2\times 2$-minors of $m_{\K}(x)$.
From now on we shall assume that $P_{\K}=\K[x]$, so that $\operatorname{rank}(m_{\K}(p))\geq 2$
at any $p\in\K^4$ and then
$V(S_{\K})=\{ p\in\K^4\ |\ \operatorname{rank}(m_{\K}(p))=2\}$.

In particular, if $p\in V(S_{\R})$ (resp. $p\in V(S_{\C})$) then $p$ is isolated in $m_{\R}^{-1}(\Sigma)$
(resp. in $m_{\C}^{-1}(\Sigma_{\C}))$. Hence, for $p\in V(S_{\R})$
the intersection index $I_p(m_{\R})$ is defined, and is equal to
the local topological degree $\deg_p(H_{\R}):(\R^4,p)\rightarrow (\R^4,\bo)$.\\[1em]

\noindent{\bf Definition.}  Set $\#\Sigma^2(m_{\R})=\sum I_p(m_{\R})$, where
$p\in V(S_{\R})=m_{\R}^{-1}(\Sigma)$.\\[1em]

Be Lemma \ref{przestawienia}, after  $\R$-linear changes of coordinates in $\C^n$ if necessary,
one may expect  $\det(A)$ not to vanish at every $p\in V(S_{\K})$, so that the ideal in $\K[x]$
generated by $S_{\K}$ and $\det(A)$ equals $\K[x]$. This justifies assumptions of the next lemma.

\begin{lemma}\label{pomocniczy}
Assume that the ideal generated by $S_{\K}$ and  $\det(A)$ equals $\K[x]$. Then
\begin{itemize}
\item[(i)] $\det(A(p))\neq 0$ at each $p\in V(S_{\K})$,
\item[(ii)] at each $p\in V(S_{\K})$, the ideal in $\Opis_{\K,p}$ generated by $S_{\K}$ equals the one
generated by $h_{44},h_{43},h_{34},h_{33}$.
\end{itemize}
\end{lemma}
\noindent{\em Proof.} Assertion $(i)$ is obvious. If $p\in V(S_{\K})$ then $\det(A(p))\neq 0$. By Lemma \ref{minory_1},
the ideal in $\Opis_{\K,p}$  generated by $S_{\K}$, i.e. by all $3\times 3$-minors of $m_{\K}(x)$, 
is in fact generated by $M_{ij}(m_{\K}(x))=h_{ij}(x)$, where $3\leq i,j\leq 4$.  $\Box$\\[1em]

Assertion \emph{(ii)} of the above lemma allows us to compute $\#\Sigma^2(m_{\R})$
by applying arguments developed in \cite[Section 3]{karolkiewiczetal}. We shall now
recall briefly the method presented there.

Let $V(S_{\C})=\{p_1, \ldots ,p_r\}$. The complex conjugation on $V(S_{\C})$ fixes $V(S_{\R})$,
so one may assume that $V(S_{\R})=\{p_1, \ldots ,p_m \}$  and
$V(S_{\C})\setminus V(S_{\R})$ is the union of pairs of conjugate points
$\{p_{m+1},\overline{p_{m+1}},\ldots, p_w,\overline{p_w}\}$, where $w=(r-m)/2$. 
Put $h_1=h_{44}, h_2=h_{43}, h_3=h_{34}, h_4=h_{33}$.

For $x=(x_1,\ldots ,x_4),\ x'=(x_1',\ldots ,x_4')$, and $1\leq
i,j\leq 4$ define
$$ T_{ij}(x,x')=\frac{h_i(x_1',\ldots,x_j, \ldots , x_4)-
   h_i(x_1',\ldots ,x_j',\ldots ,x_4)}{x_j-x_j'}.$$
It is easy to see that each $T_{ij}$ extends to a polynomial, thus
we may assume that
$T_{ij}\in\R[x,x']=\R[x_1,\ldots ,x_4,x_1',\ldots ,x_4']$.
There is the natural projection $\R[x,x']\longrightarrow \Aa_{\R}\otimes\Aa_{\R}$ given by
\[x_1^{\alpha_1}\cdots x_4^{\alpha_4}(x_1')^{\beta_1}\cdots (x_4')^{\beta_4}
\mapsto x_1^{\alpha_1}\cdots x_4^{\alpha_4}\otimes
         (x_1')^{\beta_1}\cdots (x_4')^{\beta_4}.\]
Let $T$ denote the image of $\det [T_{ij}(x,x')]$ in $\Aa_{\R}\otimes\Aa_{\R}$.

Put $d=\dim_{\R}\Aa_{\R}$. Assume that $e_1,\ldots ,e_d$ form
a basis in $\Aa_{\R}$. So $\dim_{\R} \Aa_{\R}\otimes\Aa_{\R}=d^2$ and $e_i\otimes e_j$, for $1\leq i,j\leq d$,

form a basis in $\Aa_{\R}\otimes\Aa_{\R}$. Hence there are
$t_{ij}\in\R$ such that
\[ T=\sum_{i,j=1}^d t_{ij}e_i\otimes e_j=\sum_{i=1}^d e_i\otimes\hat{e}_i,\]
where $\hat{e}_i=\sum_{j=1}^d t_{ij}e_j.$ Elements
$\hat{e}_1,\ldots ,\hat{e}_d$ form a basis in $\Aa_{\R}$. So
there are $A_1,\ldots , A_d\in\R$ such that
$1=A_1\hat e_1+\ldots  +A_d\hat e_d\mbox{ in }\Aa_{\R}$.\\[1em]

\noindent{\bf Definition.\/} For $a=a_1e_1+\ldots + a_de_d\in\Aa_{\R}$
define $\varphi_T(a)=a_1A_1+\ldots +a_d A_d.$ Hence
$\varphi_T:\Aa_{\R}\longrightarrow \R$
is a linear functional.
Let $\Phi_T$ be the bilinear form on $\Aa_{\R}$ given by
$\Phi_T(a,b)=\varphi_T(ab)$.

\begin{theorem}{ \cite[Theorem 14, p. 275]{karolkiewiczetal}}\label{wklad}
 The form $\Phi_T$ is non-degenerate and
$$\sum_{i=1}^m \deg_{p_{i}} (H_{\R})=\operatorname{signature}
\Phi_T. \ \Box$$
\end{theorem}

\begin{theorem}
Suppose that $\dim_{\R}\Aa_{\R}<\infty$ and the ideal generated by $S_{\R}$ and  $\det(A)$ equals $\R[x]$. Then
$$\#\Sigma^2(m_{\R})=\operatorname{signature}(\Phi_T).\ \Box$$
\end{theorem}
\noindent{\em Proof.} It is enough to observe that
$\#\Sigma^2(m_{\R})=   \sum_{i=1}^m I_{p_i}(m_{\R}) =        \sum_{i=1}^m\deg_{p_i}(H_{\R}).$\ $\Box$\\[1em]
{\bf Example.} Let $f=(x-2y^2+zw,y-x^2w+4z^3,zw+3w+x^2,xz+yw-4y):\R^4\rightarrow\R^4$, and $m_{\R}=Df$. 
Using {\sc Singular} \cite{greueletal}  and computer programs written by Adriana Gorzelak and Magdalena
Sarnowska - students of computer sciences at the Gda\'nsk University - one may verify that
$\dim_{\R}\Aa_{\R}=34$, other assumptions of the above theorem hold, and
$\operatorname{signature}(\Phi_T)=2$, so that $\#\Sigma^2(Df)=2$.\\[1em]
{\bf Example.} Let $f=(x-z^3,y-xzw,x^3-yz+yw,x^2+y^2+zw):\R^4\rightarrow\R^4$, and $m_{\R}=Df$. One may verify that 
$\dim_{\R} \Aa_{\R}=23$, other assumptions of the above theorem hold,
and $\operatorname{signature}(\Phi_T)=1$, so that $\#\Sigma^2(Df)=1$. Moreover $\operatorname{rank}(Df(\bo))=2$, $I_{\bo}(Df)=-1$ and
$\dim_{\R} \Aa_{\R,\bo}=3$, so that the origin is not an umbilic point.\\[1em]

If $f:\R^4\rightarrow\R^4$ is a proper polynomial mapping, there exists $R>0$ such that $|x|<R$ for each $x\in f^{-1}(\bo)$.
Denote by $\deg(f)$ the topological degree of $S_R^3\ni x\mapsto f(x)/|f(x)|\in S^3$.\\[1em]

\noindent{\bf Example.} Let $f_\pm=(x,y,z^2-w^2\pm xz+yw,\pm zw)$. Both $f_{\pm}$ are proper
and $(Df_{\pm})^{-1}(\Sigma)=f_{\pm}^{-1}(\bo)=\{\bo\}$. We have
$$Df_{\pm}=\left[ \begin{array}{cccc}1&0&0&0\\
0&1&0&0\\
\pm z&w&2z\pm x&-2w+y\\  0&0&\pm w &\pm z                   \end{array}    \right],$$
so that $Df_{\pm}(\bo)\in \Sigma\cap U$.
Then $\#\Sigma^2(Df)=I_{\bo}(Df_{\pm})=\deg_{\bo}(2z\pm x,-2w+y,\pm w,\pm z)=\mp 1$. Moreover, $\deg(f_{\pm})=\deg_{\bo}(f_{\pm})=\pm 2$.\\[1em]

\noindent{\bf Example.} Let $g_\pm=(x,y,z^2+w^2\pm xz+yw,\pm zw)$. 
Both  $g_\pm$ are proper and $(Dg_{\pm})^{-1}(\Sigma)=g_{\pm}^{-1}(\bo)=\{\bo\}$. We have
$$Dg_{\pm}=\left[ \begin{array}{cccc}1&0&0&0\\
0&1&0&0\\
\pm z&w&2z\pm x&2w+y\\  0&0&\pm w &\pm z                   \end{array}    \right],$$
so that $Dg_{\pm}(\bo)\in \Sigma\cap U$.
Then $\#\Sigma^2(Dg)=I_{\bo}(Dg_{\pm})=\deg_{\bo}(2z\pm x,2w+y,\pm w,\pm z)=\mp 1$. Moreover, $\deg(g_{\pm})=\deg_{\bo}(g_{\pm})=0$.

The above examples demonstrate  that there is only a trival linear relation of the form $A\,\#\Sigma^2(Df)+B\,\deg(f)+C=0$
for  mappings $f:\R^4\rightarrow\R^4$, contrary to the case of mappings $M\rightarrow N$, where $M$
 is closed (see \cite{stingley}).\\[1em]
{\bf Acknowledgements.} The author would like to thank Professor Osamu Saeki and Professor Kentaro Saji
for informing him about relevant results by R.~Singley and other authors.

Zbigniew SZAFRANIEC\\
Institute of Mathematics, University of Gda\'nsk\\
80-952 Gda\'nsk, Wita Stwosza 57, Poland\\
Zbigniew.Szafraniec@mat.ug.edu.pl\\

\end{document}